\documentclass[11pt]{amsart} 
\usepackage{amssymb,amsmath, amsthm} 
\usepackage{a4}

\newcommand{\Cd}{{\mathcal K}^d_0} 
\newcommand{\Lat}{{\mathcal L}^d}
\newcommand{\R}{{\mathbb R}}  
\newcommand{\N}{{\mathbb N}} 
 
\newcommand{\Z}{{\mathbb Z}}

\DeclareMathOperator{\vol}{vol} 
\DeclareMathOperator{\inter}{int}

\DeclareMathOperator{\lin}{lin}

\newtheorem{theorem}{Theorem}[section]

\newtheorem{lemma}[theorem]{Lemma}
\newtheorem{definition}[theorem]{Definition}
\newtheorem{conjecture}[theorem]{Conjecture}

\theoremstyle{definition} 
 
\numberwithin{equation}{section}

\begin{document}

\title {Successive Minima and Lattice Points}  
 
\author{Martin Henk}
\address{Martin Henk, Technische Universit\"at Wien, Abteilung f\"ur Analysis,
  Wiedner Hauptstr. 8-10/1142, A-1040 Wien, Austria} 
\email{henk@tuwien.ac.at}


\thanks{The work was completed while I  was 
        visiting the Mathematical Department of the University of Crete in
        Heraklion.  I would like to thank  
        the University of Crete for their great hospitality and support.
         }

\begin{abstract} The main purpose of this note is to prove an upper
  bound on the number of lattice points of a centrally symmetric
  convex body in terms of the successive minima of the body. This bound
  improves on former bounds and narrows the gap towards a lattice
  point analogue of Minkowski's second theorem on successive minima.   
  
  Minkowski's proof of his second theorem is rather lengthy and it was
  also criticised as obscure. We present a short proof of  Minkowski's
  second theorem on successive minima, which, however, is based on the
  ideas of Minkowski's proof.
\end{abstract}

\maketitle


\section{Introduction}
In 1896 Hermann Minkowski's fundamental and guiding book ``Geometrie der
Zahlen'' \cite{Min:geozahl} was published, which
 may be considered as the first systematic study on relations
between convex geometry, Diophantine approximation, and the theory of
quadratic forms (cf.~{\sc Gruber} \cite{Gru:geonum}). One of the
basic problems in   
geometry of numbers is to decide whether a given set in the
$d$-dimensional Euclidean space $\R^d$ contains a non-trivial lattice point
 of a $d$-dimensional lattice $\Lambda\subset \R^d$. 
With respect to the class $\Cd$ of all 0-symmetric convex bodies in 
$\R^d$ with non-empty interior  and the volume
$\vol(\cdot)$ -- $d$-dimensional Lebesgue measure -- 
Minkowski settled this problem:  
\begin{equation}
\label{eq:minkowski0}
\text{\it If $\vol(K)\geq 2^d\cdot \det\Lambda$ then  $K$
  contains a non-zero lattice point of $\Lambda$. } 
\end{equation}
Here $\det\Lambda$ denotes the determinant of the lattice $\Lambda$
and the space of all lattices $\Lambda\subset\R^d$ with
$\det\Lambda\ne 0$ is denoted by $\Lat$. 
 
Minkowski assessed  his result as ``ein Satz, der nach meinem
Daf\"urhalten zu den fruchtbarsten in der Zahlenlehre zu rechnen
ist'' (\cite{Min:geozahl}, p.~75) and indeed this theorem has many
applications (cf.~\cite{ErdGruHam:lattpoint}, sec.~3.3).   
Minkowski proved even a
stronger result, for which  we have to introduce his ``kleinstes
System von unabh\"angig gerichteten Strahlendistanzen im
Zahlengitter'' (\cite{Min:geozahl}, p.~178). 
\begin{definition} Let $K\in\Cd$ and $\Lambda\in\Lat$.
For $1\leq i\leq d$  
\begin{equation*}
 \begin{split}
 \lambda_i(K,\Lambda)=\min\big\{
    \lambda\in\R_{\geq 0}: &\,\lambda K \,\,\text{contains }\, i
    \text{ linearly independent} \\
     &\text{lattice points of }\Lambda\big\}
 \end{split}
\end{equation*}
is called the {\em $i$-th successive minimum of $K$ with respect to~$\Lambda$}. 
\end{definition}
Obviously, we have $\lambda_1(K,\Lambda)\leq
\lambda_2(K,\Lambda)\leq\cdots\leq\lambda_d(K,\Lambda)$ and 
the first successive minimum $\lambda_1(K,\Lambda)$ 
 is the smallest dilation factor such that $\lambda_1(K,\Lambda)\, K$
contains a non-zero lattice point. With this notation Minkowski's first
theorem on successive minima reads (cf.~\cite{Min:geozahl}, pp.~75)
\begin{theorem}[Minkowski] Let  $K$ in
  $\Cd$ and $\Lambda\in\Lat$. Then 
\begin{equation*}
 \lambda_1(K,\Lambda)^d \vol(K)\leq 2^d\,\det\Lambda.
\end{equation*}
\label{thm:Minkowski_first}
\end{theorem}
So  $\vol(K)\geq 2^d\det\Lambda$ implies $\lambda_1(K,\Lambda)\leq 1$,
and we get \eqref{eq:minkowski0}. 
Minkowski's second theorem on successive minima is a deep improvement of the
first one and says (cf.~\cite{Min:geozahl}, pp.~199)
\begin{theorem}[Minkowski] Let  $K \in
  \Cd$ and $\Lambda\in\Lat$. Then
\begin{equation*}     
 \lambda_1(K,\Lambda)\cdot\lambda_2(K,\Lambda)\cdot\ldots\cdot\lambda_d(K,\Lambda)\cdot \vol(K)\leq 2^d\,\det\Lambda.
\end{equation*}
\label{thm:Minkowski_second}
\end{theorem}
This inequality is best possible. For instance, with respect to the
integral lattice $\Z^d$, each box with axes
parallel to the coordinate axes gives equality. Although Theorem
\ref{thm:Minkowski_second} has not so many applications as the first
theorem on successive minima, it shows a beautiful relation
between the volume of $K$ and the expansion of $K$ with respect to
independent lattice directions of a lattice. 
The importance of Theorem
\ref{thm:Minkowski_second} is also reflected in the number of
different proofs, see e.g.~ 
{\sc  Bambah,  Woods \& Zassenhaus}  \cite{BamWooZas:succmin},  
{\sc Cassels} \cite{Cas:geonum},  
{\sc Danicic}  \cite{Dan:succmin}, 
{\sc  Davenport} \cite{Dav:succmin}, 
{\sc  Estermann}  \cite{Est:succmin},  
{\sc Siegel} \cite{Sie:geonum} and 
{\sc  Weyl}  \cite{Wey:succmin}.

In \cite{BetkeHenkWills:successive_minima} it was conjectured that
an inequality analogue  to Theorem \ref{thm:Minkowski_second} holds for the
lattice point enumerator $ \#(K\cap\Lambda)$. 
More precisely,
\begin{conjecture} Let  $K
  \in \Cd$ and $\Lambda\in\Lat$. Then
\begin{equation}
   \#(K\cap\Lambda) \leq \prod_{i=1}^d \left\lfloor
     \frac{2}{\lambda_i(K,\Lambda)}+1\right\rfloor.  
\label{eq:conj}
\end{equation}
\label{conj:bhw}
\end{conjecture}  
Here $\lfloor x \rfloor$ denotes the smallest integer not less than
$x$. An analogous statement to the first theorem of Minkowski on
successive minima was already shown in
\cite{BetkeHenkWills:successive_minima}, namely 
\begin{equation}
\#(K\cap\Lambda) \leq 
    \left\lfloor
        \frac{2}{\lambda_1(K,\Lambda)}+1\right\rfloor^d.
\label{eq:lat_min_one}
\end{equation}
It seems to be worth mentioning that if Conjecture \ref{conj:bhw} 
 were true then  we could write by the definition of the Riemann
 integral   
\begin{equation*}
 \frac{\vol(K)}{\det\Lambda}=\lim_{r\to 0} r^d\#(K\cap r\Lambda) 
 \leq \lim_{r\to 0} \prod_{i=1}^d r \left\lfloor
     \frac{2}{\lambda_i(K,r\Lambda)}+1\right\rfloor =\prod_{i=1}^d 
     \frac{2}{\lambda_i(K,\Lambda)}.
\end{equation*}
Thus Conjecture \ref{conj:bhw} implies Minkowski's second
theorem on successive minima (Theorem \ref{thm:Minkowski_second}). 
In \cite{BetkeHenkWills:successive_minima} the validity of the
conjecture was proven in the  case $d=2$. Moreover, it was shown that
an upper bound of this type 
exists, if in the above product
$\frac{2}{\lambda_i(K,\Lambda)}$ is replaced by
$\frac{2\,i}{\lambda_i(K,\Lambda)}$. So, roughly speaking,
\eqref{eq:conj} holds up to a factor $d!$. Here we shall improve this
bound.
\begin{theorem} Let $d\geq 2$,  $K\in \Cd$ and $\Lambda\in\Lat$. Then    
\begin{equation*}
   \#(K\cap\Lambda) < 2^{d-1} \prod_{i=1}^d \left\lfloor
     \frac{2}{\lambda_i(K,\Lambda)}+1\right\rfloor.  
\end{equation*}
\label{thm:main}
\end{theorem}
The proof of this theorem will be given is the next section.  
Minkowski's original proof
(\cite{Min:geozahl}, 199-218) of his second theorem on successive minima
was sometimes
criticised as lengthy and obscure (cf.~\cite{davenport:collected_mink},
p.91). One reason might be that in the scope of the proof he also proves
many basic facts about the volume of a convex body,
like the computation of the volume through successive integrations, etc., 
which cloud a little bit the simple and nice geometrical ideas of his
proof.   Based on these ideas we present a short proof of Theorem
\ref{thm:Minkowski_second} in the last section.

For more information on lattices, successive
minima and their role in the geometry of numbers we refer to the books
of {\sc Erd\"os, Gruber and Hammer} \cite{ErdGruHam:lattpoint}, 
 {\sc Gruber  and Lekkerkerker} \cite{GruLek:geonum} and the survey of
 {\sc Gruber} \cite{Gru:geonum}.  For an elementary introduction to
 the geometry of numbers see \cite{OldsLaxDavidoff:geometry_of_numbers}.

\section{Proof of Theorem \ref{thm:main}} 
Before giving the proof we list some basic facts on
lattices, for which we refer to \cite{GruLek:geonum}.
Every lattice $\Lambda\in\Lat$ can be
written as $\Lambda=A\Z^d$, where $A$ is a non-singular $(d\times
d)$-matrix, 
i.e., $A\in {\rm GL}(d,\R)$. In particular  we have 
$\lambda_i(K,\Lambda)=\lambda_i(A^{-1}K,\Z^d)$ and
$\#(K\cap\Lambda)=\#(A^{-1}K\cap\Z^d)$. A lattice
$\widetilde{\Lambda}\in\Lat$ is called a sublattice of
$\Lambda\in\Lat$ if $\widetilde{\Lambda}\subset\Lambda$. For 
$a,\overline{a}\in\Lambda$ and a sublattice
$\widetilde{\Lambda}\subset\Lambda$ we write  
$$ 
   a\equiv \overline{a}\bmod \widetilde{\Lambda} \Leftrightarrow  
  (a-\overline{a})\in\widetilde{\Lambda}. 
$$  
 In words, $a,\overline{a}$ belong to the same residue class (coset) of
$\Lambda$ with respect to~$\widetilde{\Lambda}$.  We note that there are precisely
$\det\widetilde{\Lambda}/\det\Lambda$ different residue classes of
$\Lambda$ with respect to~$\widetilde{\Lambda}$. For every set of
$d$-linearly independent lattice points 
$a^1,\dots,a^d$ of a lattice $\Lambda$ there
exists a basis $b^1,\dots,b^d$ of $\Lambda$ such that
$\lin\{a^1,\dots,a^i\}=\lin\{b^1,\dots,b^i\}$,
where $\lin$ denotes the linear hull. In particular, given $d$
linearly independent lattice vectors $z^i\in\Z^d$, $1\leq i\leq d$,
with $z^i\in\lambda_i(K,\Z^d)\,K$ then there exists an unimodular
matrix $U$, i.e., $U\in {\rm GL}(d,\R)\cap \Z^{d\times
  d}$, such that 
\begin{equation}
    Uz^i\in\left(\lambda_i(UK,\Z^d)\,UK\right)
   \cap\lin\{e^1,\dots,e^i\},\quad 1\leq
    i\leq d,
\label{eq:succ_minima}
\end{equation}
where $e^i\in\R^d$ denotes the $i$-th unit vector. 
Furthermore we note that for $d$  linearly independent lattice points
$a^1,\dots,a^d$ of a lattice
$\Lambda\in\Lat$ satisfying  $a^i\in\lambda_i(K,\Lambda)\,K$, the
definition of the successive minima implies   
\begin{equation} 
     \inter\left(\lambda_i(K,\Lambda)\,K\right)\cap \Lambda
     \subset\lin\{0,a^1,\dots,a^{i-1}\}\cap \Lambda,\quad 1\leq i\leq d,
\label{eq:succ_cons}
\end{equation} 
where $\inter$ denotes the interior.

For the proof of Theorem \ref{thm:main} we need the following simple
lemma. 
\begin{lemma} Let $K\in\Cd$, $\Lambda\in\Lat$ and let $\widetilde{\Lambda}$ be a sublattice of $\Lambda$. Then  
$$
  \#\left(K\cap\Lambda\right)\leq 
       \frac{\det\widetilde{\Lambda}}{\det\Lambda}\,\#\left(2\,K\cap\widetilde{\Lambda}\right).
$$ 
\label{lem:lat_suc} 
\end{lemma}
\begin{proof} Let $m=\#(2\,K\cap\widetilde{\Lambda})$ and
  suppose there exist at least $m+1$ different lattice points
  $a^1,\dots,a^{m+1}\in K\cap\Lambda$ such that $a^i\equiv
  a^1\bmod\widetilde{\Lambda}$, $1\leq i\leq m+1$. Then we have  
$$
  a^i-a^1  \in (K-K)\cap\widetilde{\Lambda}=2K\cap
  \widetilde{\Lambda}, \quad 1\leq i\leq m+1,   
$$
which contradicts the assumption
$\#(2\,K\cap\widetilde{\Lambda})=m$. Thus we have shown that every
residue class  of $\Lambda$
with respect to~$\widetilde{\Lambda}$ does not contain more than $m$
points  of $K\cap\Lambda$. Since there are precisely
$\det\widetilde{\Lambda}/\det\Lambda$ different residue classes, we get the
desired bound.
\end{proof}
We remark that inequality \eqref{eq:lat_min_one} is a simple
consequence of this lemma. To see this we set $n_1=\lfloor
2/\lambda_1(K,\Lambda)+1\rfloor$ and
$\widetilde{\Lambda}=n_1\Lambda$. Next we observe that $a\in
 2K\cap\widetilde{\Lambda}$ implies that $\frac{1}{n_1}a\in
 \inter(\lambda_1(K,\Lambda)\,K)\cap\Lambda$ and from \eqref{eq:succ_cons} we
 conclude $\#(2K\cap\widetilde{\Lambda})=1$. Thus Lemma
 \ref{lem:lat_suc} gives 
\begin{equation*}
    \#(K\cap\Lambda) \leq \frac{\det\widetilde{\Lambda}}{\det\Lambda}=(n_1)^d  =\left\lfloor
     \frac{2}{\lambda_1(K,\Lambda)}+1\right\rfloor^d.
\end{equation*} 
Next we come to the proof of Theorem \ref{thm:main}.
\begin{proof}[Proof of Theorem \ref{thm:main}] W.l.o.g.~let
  $\Lambda=\Z^d$ and we may assume that (cf.~\eqref{eq:succ_minima}
  and \eqref{eq:succ_cons}) 
\begin{equation}
    \inter\left(\lambda_i(K,\Z^d)K\right)\cap\Z^d
   \subset\lin\{0,e^1,\dots,e^{i-1}\}\cap\Z^d,\quad 
   1\leq i\leq d.
\label{eq:proof_assum}
\end{equation}
For abbreviation we set
$q_i=\lfloor\frac{2}{\lambda_i(K,\Lambda)}+1\rfloor$, $1\leq i\leq
d$, and first we determine  $d$ numbers $n_i\in\N$ such that 
\begin{equation}
\begin{split}
    n_d  =q_d, \quad 
    q_i \leq n_i
    < 2\,q_i, \quad  \text{ and }\quad n_{i+1} \text{ divides } n_i,\quad\,1\leq i\leq d-1. 
\end{split}
\label{eq:proof_numbers}
\end{equation}
Suppose we have already found 
  $n_d,\dots,n_{k+1}$ with these properties. In order to determine
  $n_k$ we distinguish two cases. If $n_{k+1}\geq q_k$ we 
  set $n_k=n_{k+1}$. Since $q_k\geq q_{k+1}$ we obtain $q_k\leq
  n_k=n_{k+1}<2\,q_{k+1}\leq 2\,q_k$. Otherwise, if $n_{k+1}<q_k$ let 
 $q_k=m\cdot n_{k+1}+r$ with $m\in\N$, $m\geq 1$, and $0\leq r<n_{k+1}$. In this
 case we set $n_k=q_k+n_{k+1}-r$ and obviously, $n_k$ meets the
 requirements of \eqref{eq:proof_numbers}. 

Now let $\widetilde{\Lambda}\subset \Z^d$ be the lattice generated by the vectors 
$n_1\,e^1,n_2\,e^2,\dots,n_d\,e^d$. Then we have
$\det\widetilde{\Lambda}/\det\Lambda=n_1\cdot n_2\cdot\ldots\cdot n_d$
and together with the upper bounds on the
 the numbers $n_i$,  Lemma \ref{lem:lat_suc} gives
\begin{equation}
 \#\left(K\cap\Lambda\right)\leq \#\left(2\,K\cap\widetilde{\Lambda}\right)\,\prod_{i=1}^d n_i
 < \#\left(2\,K\cap\widetilde{\Lambda}\right)\,2^{d-1}
 \prod_{i=1}^d\left\lfloor\frac{2}{\lambda_i(K,\Lambda)}+1\right\rfloor.
\label{eq:implies_thm}
\end{equation}
Hence, in order to verify the theorem,  it suffices to show $2K\cap\widetilde{\Lambda}=\{0\}$.  
Suppose there exists a $g\in 2K\cap\widetilde{\Lambda}\setminus\{0\}$
and let $k$ 
be the largest index of a non-zero coordinate of $g$, i.e., $g_k\ne 0$
and $g_{k+1}=\cdots =g_d=0$. Then we may write   
$$ 
       g=z_1\, (n_1e^1)+z_2\,(n_2e^2) + \cdots + z_k\,(n_ke^k)\in 2K  
$$ 
for some $z_i\in\Z$. Since $n_{k}$ is a divisor of
$n_1,\dots,n_{k-1}$ and since 
$2/n_k<\lambda_k(K,\Z^d)$ (cf.~\eqref{eq:proof_numbers}) we obtain   
$$ 
    \frac{1}{n_k} g \in \left(\frac{2}{n_k} K\right) \cap
      \Z^d\subset  
    \inter(\lambda_k(K,\Z^d)K)\cap\Z^d. 
$$  
However, since $g_k\ne 0$  this relation violates
 \eqref{eq:proof_assum}. Thus we have
$2K\cap\widetilde{\Lambda}=\{0\}$ and the theorem is proven.
\end{proof}

\section{Proof of Theorem \ref{thm:Minkowski_second} }
Minkowski's proof of his  second theorem on successive minima can be found
in his book ``Geometrie der Zahlen'' (\cite{Min:geozahl},
199--219) and for an English translation we refer to
\cite{hancock:geonum}, 570--603. 

\noindent
\begin{proof}[Proof of Theorem \ref{thm:Minkowski_second}\,\,{\rm
    (following Minkowski)}]
  Again w.l.o.g.~we may assume that  $\Lambda=\Z^d$. For convenience we write 
 $\lambda_i=\lambda_i(K,\Z^d)$ and set $K_i=\frac{\lambda_i}{2}K$. 
 Furthermore, we assume that $z^1,\dots,z^d$ are $d$ linearly
 independent lattice points with  $z^i\in \lambda_i K\cap \Z^d$ and
 $\lin\{z^1,\dots,z^i\}=\lin\{e^1,\dots,e^i\}$, $1\leq i\leq d$,
 (cf.~\eqref{eq:succ_minima}). For short, we denote the linear space $\lin\{e^1,\dots,e^i\}$ by $L_i$.  

  For an integer $q\in \N$ let   $M_q^d=\{z\in\Z^d : |z_i|\leq q,\,
  1\leq i\leq d\}$ 
  and for $1\leq j\leq d-1$ let $M_q^j=M_q^d\cap L_j$.   
  Since $K$ is a bounded set  there exists a constant $\gamma$, 
  only depending on $K$, such that 
\begin{equation}
        \vol(M_q^d+K_d)\leq (2q+\gamma)^d.
\label{eq:second_one}
\end{equation}
By the definition of $\lambda_1$ we have 
$(z+\inter(K_1))\cap(\overline{z}+\inter(K_1))=\emptyset$ for two
different lattice point $z,\overline{z}\in\Z^d$, because otherwise we
would get the contradiction
$z-\overline{z}\in(\inter(K_1)-\inter(K_1))\cap\Z^d=\inter(K_1-K_1)\cap\Z^d=
\inter(\lambda_1\,K)\cap\Z^d=\{0\}$. Thus we have  
\begin{equation}
 \vol(M_q^d+K_1)= (2q+1)^d\vol(K_1) = (2q+1)^d\left(\frac{\lambda_1}{2}\right)^d \vol(K).
\label{eq:second_two}
\end{equation}
In the following we shall show that for $1\leq i\leq d-1$
\begin{equation}
 \vol(M_q^d+K_{i+1})\geq \left(\frac{\lambda_{i+1}}{\lambda_i}\right)^{d-i}
                        \vol(M_q^d+K_i).
\label{eq:second_three}
\end{equation}
To this end we may assume  $\lambda_{i+1}>\lambda_i$ and let
$z,\overline{z}\in\Z^d$, which differ in the last $d-i$ coordinates, i.e., $(z_{i+1},\dots,z_{d})\ne (\overline{z}_{i+1},\dots,\overline{z}_{d})$. Then 
\begin{equation}
        \left[z+\inter(K_{i+1})\right]\cap \left[\overline{z}+\inter(K_{i+1})\right]=\emptyset.
\label{eq:cut}
\end{equation}
Otherwise the $i+1$ linearly independent lattice points $z-\overline{z},z^1,\dots,z^i$ belong to the interior of $\lambda_{i+1}K$ which contradicts the minimality of $\lambda_{i+1}$. Hence we obtain from \eqref{eq:cut}
\begin{equation*}
\begin{split}
     \vol\left(M_q^d+K_{i+1}\right) & =
     (2q+1)^{d-i}\,\vol\left(M_q^i+K_{i+1}\right), 
      \\ 
    \vol\left(M_q^d+K_{i}\right) & =(2q+1)^{d-i}\,
     \vol\left(M_q^i+K_{i}\right). 
\end{split}
\end{equation*}
and in order to verify \eqref{eq:second_three} it suffices to show
\begin{equation}
   \vol\left(M_q^i+K_{i+1}\right)\geq
   \left(\frac{\lambda_{i+1}}{\lambda_i}\right)^{d-i} \vol(M_q^i+K_i).
 \label{eq:second_five}
\end{equation} 
Let $f_1,f_2:\R^d\to \R^d$ be the linear maps given by  
\begin{eqnarray*}
  f_1(x)&=&\left( \frac{\lambda_{i+1}}{\lambda_i} x_1,\dots,
                \frac{\lambda_{i+1}}{\lambda_i} x_i,  x_{i+1},\dots,x_d\right)^\intercal,\\
 f_2(x)&=&\left(x_1,\dots,x_i,\frac{\lambda_{i+1}}{\lambda_i} x_{i+1},\dots,
                \frac{\lambda_{i+1}}{\lambda_i} x_d\right)^\intercal.
\end{eqnarray*}
Since $M^i_q+K_{i+1}=f_2(M_q^i+f_1(K_i))$ we get 
$$
  \vol\left(M_q^i+K_{i+1}\right)=
   \left(\frac{\lambda_{i+1}}{\lambda_i}\right)^{d-i} \vol(M_q^i+f_1(K_i)) 
$$
and for the proof of \eqref{eq:second_five} we have to show 
\begin{equation}
 \vol\left(M_q^i+f_1(K_i)\right) \geq \vol\left(M_q^i+K_i\right).
\label{eq:second_six}
\end{equation}
To this end  let $L_i^\perp$ be the $(d-i)$-dimensional
orthogonal complement of $L_i$. Then it is easy to see that   
for every $x\in L_i^\perp$ there exists a $t(x)\in L_i$ with 
$K_i\cap (x+L_i)\subset (f_1(K_i)\cap (x+L_i))+t(x)$ and so
$$
    \left(M^i_q+K_i\right)\cap\left(x+L_i\right)\subset 
\left[\left(M^i_q+f_1(K_i)\right)\cap\left(x+L_i\right)\right]+t(x).
$$ 
Thus we get  
\begin{eqnarray*}
               \vol(M_q^i+K_i)&=&\int_{x\in L_i^{\perp}} 
      \vol_i\left((M_q^i+K_i)\cap(x+L_i)\right){\rm d}\,x \\ 
                           &\leq &
                             \int_{x\in L_i^{\perp}} 
       \vol_i\left((M_q^i+f_1(K_i))\cap(x+L_i)\right){\rm d}\,x \\
                & = &  
                             \vol(M_q^i+f_1(K_i)), 
\end{eqnarray*}
where $\vol_i(\cdot)$ denotes the $i$-dimensional volume. This shows
\eqref{eq:second_six} and so we have verified \eqref{eq:second_three}.
Finally, it follows from \eqref{eq:second_one}, \eqref{eq:second_two} and \eqref{eq:second_three}
\begin{equation*}
 \begin{split}
 (2q+\gamma)^d  & \geq 
               \vol\left(M^d_q+K_d\right) \geq 
               \left(\frac{\lambda_d}{\lambda_{d-1}}\right) 
    \vol\left(M^d_q+K_{d-1}\right) \\ 
    & \geq \left(\frac{\lambda_d}{\lambda_{d-1}}\right)\left(\frac{\lambda_{d-1}}{\lambda_{d-2}}\right)^2 \vol\left(M^d_q+K_{d-2}\right)
   \geq \cdots\cdots \\ & \geq  
    \left(\frac{\lambda_d}{\lambda_{d-1}}\right)\cdot 
    \left(\frac{\lambda_{d-1}}{\lambda_{d-2}}\right)^2\cdot
    \ldots\cdot 
    \left(\frac{\lambda_{2}}{\lambda_{1}}\right)^{d-1}\,
    \vol\left(M^d_q+K_1\right) \\[0.5ex]
  & = \lambda_d\cdot\ldots\cdot\lambda_1\,
  \frac{\vol(K)}{2^d}\,(2q+1)^d 
 \end{split}
\end{equation*}
and so 
$$
   \lambda_1\cdot\ldots\cdot\lambda_d\,
   \vol(K)
 \leq 2^d \cdot\left(\frac{2q+\gamma}{2q+1}\right)^d.
$$
Since this holds for all $q\in\N$ the theorem  is proven.
\end{proof}

\bibliographystyle{amsalpha} 
\bibliography{martins_bibliography}

\providecommand{\bysame}{\leavevmode\hbox to3em{\hrulefill}\thinspace}
\begin{thebibliography}{BHW93}

\bibitem[BHW93]{BetkeHenkWills:successive_minima}
U.~Betke, M.~Henk, and J.M. Wills, \emph{Successive-minima-type inequalities},
  Discrete Comput. Geom. \textbf{9} (1993), 165 -- 175.

\bibitem[BWZ65]{BamWooZas:succmin}
R.~P. Bambah, A.~C. Woods, and H.~Zassenhaus, \emph{Three proofs of
  {M}inkowski's second inequality in the {G}eometry of {N}umbers},
  J.~Austral.~Math.~Soc. \textbf{5} (1965), 453--462.

\bibitem[Cas59]{Cas:geonum}
J.~W.~S. Cassels, \emph{An {I}ntroduction to the {G}eometry of {N}umbers},
  Springer, Berlin, 1959.

\bibitem[Dan69]{Dan:succmin}
I.~Danicic, \emph{An elementary proof of {M}inkowski's second inequality},
  J.~Austral.~Math.~Soc. \textbf{10} (1969), 177--181.

\bibitem[Dav39]{Dav:succmin}
H.~Davenport, \emph{Minkowski's inequality for the minima associated with a
  convex body}, Quarterly J.~Math. \textbf{10} (1939), 119--121.

\bibitem[Dav77]{davenport:collected_mink}
H.~Davenport, \emph{The collected works of {H}arald {D}avenport}, vol.~I,
  Academic Press, 1977.

\bibitem[EGH89]{ErdGruHam:lattpoint}
P.~Erd{\"o}s, P.M. Gruber, and J.~Hammer, \emph{Lattice points}, Longman
  Scientific \& Technical, Harlow, Essex/Wiley, New York, 1989.

\bibitem[Est46]{Est:succmin}
T.~Estermann, \emph{Note on a theorem of {M}inkowski}, J.~London Math.~Soc.
  \textbf{21} (1946), 179--182.

\bibitem[GL87]{GruLek:geonum}
P.M. Gruber and C.G. Lekkerkerker, \emph{Geometry of {N}umbers}, 2nd ed.,
  North-Holland, Amsterdam, 1987.

\bibitem[Gru93]{Gru:geonum}
P.M. Gruber, \emph{Geometry of numbers}, Handbook of {C}onvex {G}eometry (P.M.
  Gruber and J.M. Wills, eds.), vol.~B, North-Holland, Amsterdam, 1993.

\bibitem[Han64]{hancock:geonum}
H.~Hancock, \emph{Development of the geometry of numbers}, vol.~2, Dover
  Publications, Inc., New York, 1964.

\bibitem[Min96]{Min:geozahl}
H.~Minkowski, \emph{Geometrie der {Z}ahlen}, Teubner, Leipzig-Berlin, 1896,
  Reprinted: Johnson, New York, 1968.

\bibitem[OLD00]{OldsLaxDavidoff:geometry_of_numbers}
C.D. Olds, A.~Lax, and G.~Davidoff, \emph{The {G}eometry of {N}umbers}, The
  Anneli Lax New Mathematical Library, vol.~41, The Mathematical Association of
  America, 2000.

\bibitem[Sie89]{Sie:geonum}
C.L. Siegel, \emph{Lectures on the {G}eometry of {N}umbers}, Springer,
  New-York, 1989.

\bibitem[Wey42]{Wey:succmin}
H.~Weyl, \emph{On {G}eometry of {N}umbers}, Proc.~London Math.~Soc.~(2)
  \textbf{47} (1942), 268--289.

\end{thebibliography}

\end{document}